\documentclass[10pt]{article}

\usepackage{amssymb,color}
\usepackage[tbtags]{amsmath}

\usepackage{color}

\definecolor{c20}{rgb}{0.,0.7,0.}
\definecolor{c30}{rgb}{0.,0.,1.}
\definecolor{c40}{rgb}{1,0.1,0.7}
\definecolor{c50}{rgb}{1,0,0}

\def\cE#1{\textcolor{c20}{#1}}
\def\cE#1{#1}

\def\cE#1{#1}

\newcommand{\tb}[1]{{\textcolor{blue}{#1}}}
\def\tb#1{#1}
\newcommand{\toprob}{ \stackrel{p}{\to}}

\newtheorem{theo}{Theorem}[section]
\newtheorem{sat}[theo]{Proposition}
\newtheorem{de}[theo]{Definition}
\newtheorem{lem}[theo]{Lemma}
\newtheorem{exxa}[theo]{Example}

\newtheorem{korr}[theo]{Corollary}
\newtheorem{remark}[theo]{Remark}
\newtheorem{remarks}[theo]{Remarks}

\newcommand{\netheo}[1]{{Theorem \ref{#1}}}

\newcommand{\prooftheo}[1]{ \textsc{Proof of Theorem} \ref{#1} }
\newcommand{\proofprop}[1]{\textsc{Proof of Proposition} \ref{#1}}

\newcommand{\BQN}{\begin{eqnarray}}
\newcommand{\EQN}{\end{eqnarray}}
\newcommand{\BQNY}{\begin{eqnarray*}}
\newcommand{\EQNY}{\end{eqnarray*}}

\newcommand{\BS}{\begin{sat}}
\newcommand{\ES}{\end{sat}}
\newcommand{\BT}{\begin{theo}}
\newcommand{\ET}{\end{theo}}
\newcommand{\BK}{\begin{korr}}
\newcommand{\EK}{\end{korr}}

\newcommand{\BD}{\begin{de}}
\newcommand{\ED}{\end{de}}

\newcommand{\BIT}{\begin{itemize}}
\newcommand{\EIT}{\end{itemize}}
\newcommand{\BDI}{\begin{description}}
\newcommand{\EDI}{\end{description}}

\newcommand{\BRM}{\begin{remarks}}
\newcommand{\ERM}{\end{remarks}}

\newcommand{\QED}{\hfill $\Box$}

\newcommand{\IF}{\infty}
\newcommand{\BTH}{\begin{theo}}
\newcommand{\ETH}{\end{theo}}
\newcommand{\BPR}{\begin{sat}}
\newcommand{\EPR}{\end{sat}}

\newcommand{\BEX}{\begin{exxa}}
\newcommand{\EEX}{\end{exxa}}

\newcommand{\BC}{\begin{cases}}
\newcommand{\EC}{\end{cases}}
\newcommand{\COM}[1]{}
\newcommand{\BL}{\begin{lem}}
\newcommand{\EL}{\end{lem}}

\newcommand{\kb}[1]{\boldsymbol{#1}}
\newcommand{\vk}[1]{\kb{#1}}

\newcommand{\abs}[1]{\lvert #1 \rvert}
\newcommand{\Abs}[1]{ \Biggl \lvert #1 \Biggr \rvert}

\newcommand{\E}[1]{\mbox{\rm$\vk{E}$}\{#1\}}

\newcommand{\pb}[1]{\mathbb{P}\left\{#1 \right \}}
\newcommand{\pk}[1]{\mathbb{P}\left\{#1 \right \}}
\newcommand{\R}{\mathbb{R}}

\newcommand{\inr}{\in \R}

\newcommand{\ldot}{,\ldots,}

\newcommand{\limit}[1]{\lim_{#1 \to   \infty}}

\newcommand{\todis}{\stackrel{d}{\to}}
\newcommand{\equaldis}{\stackrel{d}{=}}
\def\zmv{ \zeta_{v}^*}

\def\zmvv{Z_{v}^*}

\def\tHH{\widehat{\varrho}}

\def\polhk#1{\setbox0=\hbox{#1}{\ooalign{\hidewidth
    \lower1.5ex\hbox{`}\hidewidth\crcr\unhbox0}}}

\begin{document}
\quad\\
\quad

\centerline{\Large \bf Gaussian Approximation of Perturbed Chi-Square Risks}

\bigskip

\centerline{ \bf Krzysztof D\polhk{e}bicki\footnote{ Mathematical Institute, University of Wroc\l aw, pl. Grunwaldzki 2/4, 50-384
Wroc\l aw, Poland}
 Enkelejd  Hashorva\footnote{Department of Actuarial Science, University of Lausanne, UNIL-Dorigny 1015 Lausanne, Switzerland}  and Lanpeng Ji$^2$}

\bigskip

{\bf Abstract}: In this paper we show that the conditional distribution of perturbed chi-quare risks can be approximated by certain distributions including the Gaussian ones. Our results are of interest for conditional extreme value models and multivariate extremes as shown in three applications.



\section{Introduction}
%1. Add Sebastian Bernoulli
Let $(X_{i1},X_{i2}), i\ge 1$ be independent bivariate Gaussian random vectors with $N(0,1)$ distributed marginals and correlation coefficient  $\rho \in (-1,1)$. Clearly, they have the following  stochastic representation
\BQN
\label{eJ}
(X_{i1},X_{i2}) \equaldis (X_{i1}, \rho  X_{i1}+ \sqrt{1- \rho^2} W_i),\ \ i\ge1,
\EQN
with $W_i,i\ge 1$ being independent $N(0,1)$ random variables which are further independent of $X_{i,1}, i\ge 1$.
For fixed $m\ge 2$ we define a
 bivariate chi-square random vector $(  \zeta_1,   \zeta_2)$ by
\BQN \label{chi12}
   \zeta_1= \sum_{i=1}^m X_{i1}^2, \quad %l=1,2
   \zeta_2= \sum_{i=1}^m X_{i2}^2.
\EQN
\tb{Obviously, apart from the
case $\rho=0$, it has dependent components. By a direct analytic proof (see Appendix)
it follows that, as $v\to \IF$, the conditional risk 
\BQNY%\label{spec}
\zmv:=\frac{  \zeta_2  - \rho^2v }{2 \rho \sqrt{(1 - \rho^2)v}} \Bigl \lvert (  \zeta_1= v)
%, \quad \text{ and   }\zzmv:=\frac{ \zeta_2  - \rho v }{ \tHH }
%\Bigl \lvert ( \zeta_1 = v)
%, \quad v>0%, \quad \text{  with } \tHH:= \sqrt{1 - \rho^2}
\EQNY
can be approximated by a standard Gaussian random variable $W$,
in such a way that} %namely %for any $m\ge 1$ and  any $x\inr$
\BQN\label{eq:A}
\limit{v} \sup_{x\inr} \Abs{ \pb{ \zmv \le x }  %\limit{v} \pb{ \zzmv \le x }
- \pk{W \le x}}=0.
\EQN
Instead of conditioning on $\{  \zeta_1 =v\}$ \tb{one} can also condition on the
event $\{   \zeta_1 > v\}$. \tb{Moreover, the} same Gaussian approximation
 of $ \zeta_2$ given that $\{   \zeta_1 > v\}$ can be obtained (see \netheo{Th2} below).

\tb{The motivation of analyzing the distributional properties of
the conditional models stems both from theory- and applied-oriented problems.}
Commonly in finance and risk management applications there are few observations of risks being large. Therefore, a conditional model, which can be reasonably approximated by some known distribution functions, is valuable for statistical models; see e.g., \cite{MR2543586, balkema2007high, eng2012a, MR2739357, MR2891307, HashExt12, MR3003975, heffernan:resnick:2007, heffernan:tawn:2004, pre06153156} for various results.

Conditional limit results are also crucial  for the investigation of the asymptotic behaviour
of maximum of stationary random processes and that of maxima of
triangular arrays; see e.g.,
\cite{MR1043939, Berman82, Berman92, hue1989, HB, Pit96, MR2733939}
and references therein.
{\tb{Namely,} let $\{Y(t),t\ge 0\}$ be a stationary process with a.s.
continuous sample paths, and set $M(T):=\max_{t\in [0,T]} Y(t), T>0$.
In view of the works of Berman, a crucial condition for the study of $M(T)$
is the convergence of \tb{appropriate} conditional distributions.
For illustration purposes \tb{let us} assume that
$Y(0)$ has a distribution function in the Gumbel max-domain of attraction, i.e.,
\BQN\label{Gmax}
\limit{v} \frac{\pk{ Y(0) > v+  {x}/{w(v)}}}{\pk{ Y(0)> v}}=\exp(-x), \quad x\inr,
\EQN
with $w(\cdot)$ some positive scaling function.
\tb{
Then %the most crucial %condition %for the study of the tail behaviour of $M(T)$
the weak limit of
$ \vk{Y}_v(a,k):= w(v)\Bigl( Y( q(v) a)- v \ldot  Y( q(v) a k)- v\Bigr) \lvert (Y(0)> v)$
as
$v\to \IF$,
where $a>0$ and $q(\cdot)$ is a function satisfying
$\lim_{v\to\IF}q(v)=0$,
is crucial for the analysis of the asymptotics of the tail distribution
of $M(T)$.
For example, it applies to $Y$ being a stationary chi-square process with $m$ degrees of freedom.}
It is known then that $w(x)\equiv1/2$ and we can define $q(\cdot)$ and show the weak convergence of $\vk{Y}_v(a,k)$ if the underlying Gaussian process has a covariance function being regularly varying at 0;
 see e.g., Theorem 10.1 in \cite{Berman82}. 
A similar result for stationary chi-processes was also shown in \cite{MR1043939}.\\

\tb{Other important applications of approximations of the
conditional distributions of chi-square risks can be found in \cite{HashKab}.}
Consider $\zeta_{1,v}$ and $  \zeta_{2,v}$ to be realizations of some stationary chi-square process $ \{\zeta(t), t\ge0\}$ at threshold dependent times $t_1(v), t_2(v)$. 
 In this case we have threshold dependent correlation coefficient $\rho_v$ instead of constant $\rho$. 
 Clearly, in order to get results as those of Berman (see also \cite{MR1043939}) we need to assume that % control the speed that
 $\rho_v$ tends to 1 at certain speed.
This case has been considered in the context of maxima of chi-square triangular arrays in Theorem 2.1 in \cite{HashKab} which shows that $(  \zeta_{2,v} -  \zeta_{1,v})\lvert (  \zeta_{1,v} >v)$ can be approximated as $v \to \IF$  by a Gaussian random variable
with distribution $N(- \lambda, 4 \lambda)$, provided that
\BQN \label{bnn}
\limit{v} 2 v (1- \rho_v) = \lambda\in (0,\IF).
\EQN

Given the importance of approximation of conditional distributions for both cases that $\rho$ is a constant and $\rho=\rho_v$ changes with the threshold $v$,
in this paper we shall investigate approximations of multivariate conditional perturbed chi-square risks (see Section 2 for the definition) using ideas and techniques from extreme value theory.
Our findings provide a concrete framework for the conditional extreme value model developed in \cite{heffernan:tawn:2004, MR2775870},
and therefore 
statistical inference can be done using the conditional extreme value methodology  therein.
Since our approach is asymptotic in nature, distributional assumptions can be dropped.
\tb{This} makes the model more appealing for applications. More precisely, we shall drop any distributional assumption on $X_{i1}, i\ge1$. The Gaussianity of the components $W_i, i\ge 1$ in \eqref{eJ} seems to be crucial;
however there are specific models (see Section 3) where this assumption \tb{is} relaxed.

 In this paper we present three applications: The first one establishes the so-called Berman's sojourn limit theorem and the tail asymptotics of supremum for a class of time-changed stationary chi-square processes.
 The second one strengthens the convergence in distribution of  maxima of chi-square triangular arrays (see \cite{HashKab} and \cite{RevStat14})
to convergence of the corresponding probability density functions (pdf's). We conclude Section 4 with the third one concerning extremal behaviour of aggregated log-chi risks.

Organization of the paper: We \tb{begin} with the description of two main dependent perturbed chi-square models for our multivariate framework and then
derive conditional limit theorems for the models both with fixed parameters and with parameters that depend on the
threshold; \tb{see Section 2}.
Section 3 is \tb{devoted} to discussions.
The aforementioned applications are displayed in Section 4. Proofs of all results are relegated to Section 5 followed by a short Appendix.

\section{Model Description and Main Results}
\tb{We first introduce  the multidimensional perturbed chi-square random vectors.}
% analogous to the 2-dimensional one in \eqref{eJ}.
Let $(X_{i1} \ldot X_{i (k+1)}), 1\le i\le m$ be  $(k+1)$-dimensional random vectors with stochastic representations
\BQN
\label{eKJ}
(X_{i1}\ldot X_{i(k+1)}) \equaldis \Bigl(X_{i1}, \rho_{1}  X_{i1} +W_{i1} \ldot
\rho_{k}  X_{i1}+ W_{ik}\Bigr),\ \ \ 1\le i\le m,
\EQN
where $\rho_j\inr\setminus\{0\}, 1\le j\le k$, and $\mathcal{W}:=\{W_{ij}\}_{1\le i\le m, 1\le j\le k}$ is an $m\times k$ matrix of centered (non-standard) Gaussian random variables.
Define the {\it $k+1$-dimensional perturbed chi-square risk} $\vk{ \zeta}:=(\zeta_1 \ldot \zeta_{k+1})$ by
\BQN\label{eq:zetak}
  \zeta_1= \sum_{i=1}^m X_{i1}^2\ldot  \zeta_{k+1}= \sum_{i=1}^m X_{i (k+1)}^2.
  \EQN
In the sequel we shall consider the following framework:

{\bf Assumption A}: \tb{Random vector} $(X_{11} \ldot X_{m1})$ and the Gaussian
random matrix $\mathcal{W}$ are \tb{mutually} independent. Further, we assume
that the rows of  $\mathcal{W}$ are independent and have the same distribution as the centered $k$-dimensional Gaussian random vector $\vk{W}=(W_1 \ldot W_k)$. Suppose that $  \zeta_1$ has a continuous distribution function $G$ which has support on $[0,\IF)$.\\

Note  that we do not assume  $X_{11} \ldot X_{m1}$ to be independent or normally distributed. If they %$X_{11} \ldot X_{m1}$
 are independent and normally distributed with variance 1, and further for any $1\le i\le m$, $W_{i,j}$ has variance $1-\rho_j^2\in(0,1)$ for all $1\le j\le k$, then $\vk{ \zeta}$ is  the {\it (classical) chi-square risk}. %Further, the assumption that $G$ is continuous can be removed, however for notational simplicity we shall keep it.

In order to obtain an approximation for the {\it conditional perturbed chi-square risk} $(\zeta_2 \ldot \zeta_{k+1})| ( \zeta_1> v)$
we need to impose an asymptotic tail condition on $  \zeta_1$. We shall assume that $  \zeta_1$ has distribution function $G$ in the Gumbel max-domain of attraction with positive scaling function $w(\cdot)$, i.e.,
\BQN
\limit{v} \frac{1-G( v+  {x}/{w(v)} )}{1-G( v)}=\exp(-x), \quad x\inr.
\EQN
\tb{We refer to} \cite{Res1987, Faletal2010} for more details on Gumbel max-domain of attraction.
Due to the restrictions imposed by our dependence structure, not every possible scaling function $w(\cdot)$ can be considered. Thus we assume that
\BQN\label{w}
\limit{v} (\sqrt{v} w(v))^{-1}=2c\in [0, \IF).
\EQN
 Next, we state  our first result which shows convergence in distribution of the conditional perturbed chi-square risk. In what follows, the standard notation $\todis$ and $\toprob$
  denote convergence in distribution and  convergence in probability, respectively, when the argument tends to infinity.

\BT \label{Th2}
Let $\vk{ \zeta}:=(\zeta_1 \ldot \zeta_{k+1})$ be a perturbed chi-square risk given
in \eqref{eq:zetak}. Assume that {Assumption A} is satisfied,
and let $\vk{U}=(U_1 \ldot U_k)$ \tb{has} the same distribution as the centered Gaussian random vector $\vk{W}$. Then %If $\rho_i\not=0, i\le m$, then  for all $m\ge 2$ %the convergence in distribution
\BQN\label{eq:Th1U1}
\Biggl( \frac{  \zeta_2 - \rho_1^2 v}{2 \rho_1 \sqrt{v}} \ldot \frac{ \zeta_{k+1} - \rho_k^2 v}{2 \rho_k \sqrt{v}} \Biggr)  \Bigl \lvert (  \zeta_1=v)  \equaldis \vk{\zeta}_v^*
\todis  \vk{U}
 \EQN
holds, % as $v\to \IF$
 where $\vk{\zeta}_v^*, v>0$ are defined on the same probability space as $\vk{ \zeta}$. Further,
%\eqref{eq:Th1U1} holds.\\  %with $\vk{U}\equaldis \vk{W}$.\\
if $G$ is in the Gumbel max-domain of attraction with some positive scaling function $w(\cdot)$ which satisfies \eqref{w}, then
%we have %the following joint convergence in distribution
\BQN
\Biggl( w(v)(  \zeta_1 - v),  \frac{  \zeta_2 - \rho_1^2 v}{2 \rho_1 \sqrt{v}} \ldot \frac{ \zeta_{k+1} - \rho_k^2 v}{2 \rho_k \sqrt{v}} \Biggr)  \Bigl \lvert (  \zeta_1>v) \equaldis  \Bigl(\widetilde{\zeta}_v, \widetilde{\vk{\zeta}}_v\Bigr) \todis (E, \rho_1 c E+U_1 \ldot
 \rho_k c E+U_k),
\EQN
%as $v\to \IF,$
where $E$ is a unit exponential random variable independent of $\vk{U}$ and $\Bigl(\widetilde{\zeta}_v, \widetilde{\vk{\zeta}}_v\Bigr), v>0$ are defined on the same probability space as $\vk{ \zeta}$.
\ET

\begin{remark}
 In view of  \netheo{Th2}, Proposition 4.1 in \cite{MR2797990} implies that the random vector $\vk{ \zeta}$ has asymptotically independent components, i.e.,
\BQN
\limit{v} \pk{ \zeta_j> v  \lvert  \zeta_i>v}=0
\EQN
for any pair $(i,j)$ of different indices; see \cite{RevStat14} for a similar result.
\end{remark}

Our second result is concerned with the {\it threshold dependent perturbed $k+1$-dimensional chi-square risk} $\vk{ \zeta}_v:=(\zeta_{1,v} \ldot \zeta_{k+1,v})$ which is defined similarly as \eqref{eq:zetak} with $\rho_{j,v}, v>0$ instead of $\rho_j$ and Gaussian random matrices $\mathcal{W}_v, v>0$ instead of
  $\mathcal{W}$. Note that $\zeta_{1,v}=\zeta_1$. 
For $\rho_{j,v}$'s we shall impose the following conditions \tb{(compare with \eqref{bnn})}: 
\BQN\label{vww}
\limit{v} 4 v w(v) (1 - \rho_{j,v}) = \lambda_j \in [0,\IF), \quad 1\le j\le k.
\EQN
Note in passing that $\limit{v} v w(v)=\IF$, hence the above condition implies that $\limit{v} \rho_{j,v}=1$. Also note that if $G$ is a chi-square distribution, then $w(v)\equiv1/2$ and  thus \eqref{vww} reduces to \eqref{bnn}.

\BT \label{Th2B}
Let  $\vk{ \zeta}_v:=(\zeta_{1,v} \ldot \zeta_{k+1,v}), v>0$ be a family of threshold dependent perturbed chi-square risks with correlation coefficients $\rho_{j,v}\in\R/\{0\}, 1\le j\le k, v>0.$
Denote the first row of $\mathcal{W}_v$ by $\vk{W}_{v}$ and assume that Assumption A holds for every $v>0$.
\tb{Let $G$ be} in the Gumbel max-domain of attraction with scaling function $w(\cdot)$. Then

 i)  Assume that condition \eqref{vww} is satisfied \tb{and} % If further
\BQN \label{eqProb}
w(v)\sqrt{v  } \vk{W}_{v} \todis \vk{U}, \quad   \sqrt{w(v)} \vk{W}_{v} \toprob \vk{0}=(0 \ldot 0)\inr^k %, \quad v\to \IF
\EQN
holds for a random vector $ \vk{U}\in\R^k$. \tb{Then for any $x\inr$}
\BQN\label{eq:Th2A}
\Biggl(  w(v)(\zeta_{2,v}- v) \ldot   w(v)(\zeta_{k+1, v}- v)\Biggr) \Bigl \lvert   \left(\zeta_{1}= v+ \frac{x}{w(v)}\right) \todis \left( 2 U_1  - \frac{\lambda_1}{2}+x  \ldot 2 U_k - \frac{\lambda_k}{2} +x\right).%=: \vk{U}_{\vk{\lambda}}+(x \ldot x)%, \quad v\to \IF
\EQN
ii) If \eqref{eq:Th2A} holds for any $x \in [0,\IF)$, then
\BQN\label{eq:Th2B}
\Biggl(  w(v)(\zeta_{1}- v), w(v)(\zeta_{2,v}- v) \ldot   w(v)(\zeta_{k+1, v}- v)\Biggr) \Bigl \lvert   (\zeta_{1}> v) \todis
%(E, E+ \vk{U}_{\vk{\lambda}}),%
\left(E,E+2 U_1  - \frac{\lambda_1}{2} \ldot E+2 U_k - \frac{\lambda_k}{2} \right),
\EQN
with $E$ being a unit exponential random variable independent of $\vk{U}$.
% and under Model B we set $\vk{U}:=(U_1 \ldot U_k)=
%(\mathcal{R}_1 O_{11} \ldot \mathcal{R}_k O_{k1})$.

\ET

An immediate consequence of the above result is the following interesting limit relationship.
%, see \cite{eng2012a} for a more general framework.

\BK
Under the assumptions and notation of ii) in \netheo{Th2B} we have
\BQN\label{eq:Th2:2b}
&&\limit{v} \sup_{(x_1 \ldot x_{k}) \inr^k} \Abs{ \pb{  w(v)(\zeta_{2,v} -   \zeta_{1})\le x_1 \ldot
  w(v)(\zeta_{k+1,v} -  \zeta_{1}) \le x_{k} \Bigl \lvert  \zeta_{1} > v} \nonumber\\
  &&\ \ \ \ \ \ \ \qquad\qquad\qquad- \pb{2U_1-\frac{\lambda_1}{2} \le x_1 \ldot 2U_{k}-\frac{\lambda_k}{2}\le x_{k}}}=0.
\EQN
\EK
 The claim in \eqref{eq:Th2:2b} is of interest for statistical modeling; results in this direction are already available for some other interesting models  (see \cite{eng2012a}).

\begin{remarks}\label{RemTh2} a) The relation between \eqref{eq:Th2A} and \eqref{eq:Th2B} is known from several works of Berman; see e.g., \cite{Berman82} where additional conditions on the scaling function $w(\cdot)$ are  imposed.\\
b) Assume $\vk{\zeta}_v=(\zeta_{1,v},\zeta_{2,v}), v>0$ to be a family of 2-dimensional threshold dependent chi-square risks with $Var(W_{i1, v}) =1 - \rho_{1,v}^2\in(0,1), 1\le i\le m$. Clearly $w(x)=1/2$. Then from \eqref{vww} we have that  \eqref{eqProb} holds with $U_1= \sqrt{\lambda_1} V_1/2 $.\\
c) The proof of \eqref{eq:A} shows that under the assumptions of b), similar convergence as in \eqref{eq:Th2A} also holds %not only in distribution, but
for the corresponding pdf's.
\end{remarks}

\section{Discussions}

As we can see from \eqref{babo} in the proof of \netheo{Th2} that the symmetry property of Gaussian random variables plays a crucial role. In this section, we are mainly concerned with two tractable models relaxing the Gaussian assumptions.\\

First, we consider a bivariate perturbed chi-square risk  $(\zeta_1,\zeta_2)$ as in \eqref{eq:zetak}. We drop the Gaussian assumption on $W_{i,1}, 1\le i\le m$ in \eqref{eKJ} and assume  that $( X_{11} \ldot X_{m1})$ is a random vector with polar representation
$$ ( X_{11} \ldot X_{m1})= R (O_1 \ldot O_m),$$
where $R>0$ is a random variable with infinite upper endpoint, and $(O_1 \ldot O_m)$ is an independent of $R$  random vector such that
$\sum_{i=1}^2 O_i^2=1$ almost surely.
Since
\BQNY
\zeta_2&= &\rho_1^2\sum_{i=1}^m X_{i1}^2 + 2 \rho_1  \sum_{i=1}^m X_{i1} W_{i,1} +
\tb{\sum_{i=1}^m W_{i,1}^2}\\
&=& \rho_1^2R^2+ 2 R \rho_1   \sum_{i=1}^m O_i W_{i,1} +  \sum_{i=1}^m W_{i,1}^2
\EQNY
we obtain
\BQNY
\frac{\zeta_2 -\rho_1^2v}{2\rho_1\sqrt{ v}} \bigl \lvert (\zeta_1=v) \equaldis
\sum_{i=1}^m O_i W_{i,1} + \frac{  \sum_{i=1}^m W_{i,1}^2}{2\rho_1 \sqrt{v}}
\todis \sum_{i=1}^m O_i W_{i,1}, \quad v\to \IF.
\EQNY
If further $\vk{O}= (O_1 \ldot O_m)$ is uniformly distributed on the unit sphere of $\R^m$, % and $W_{i,1},i\le m$ are independent with distribution function $F$,
then $ \sum_{i=1}^m O_i W_{i,1} \equaldis O_1\sqrt{ \sum_{i=1}^m  W_{i,1}^2}$ which is in general not Gaussian. Clearly, it is Gaussian distributed if $X_{i,1}, 1\le i\le m$ are i.i.d. standard Gaussian random variables and $W_{i,1}, 1\le i\le m$ are i.i.d. Gaussian random variables with variance $1-\rho_1^2$.

 Another tractable model is obtained from \eqref{eq:zetak} by %dropping assumptions on $(X_{11} \ldot X_{m1})$
 restricting the following conditions on the random matrix $\mathcal{W}$.
Suppose that each column  $(W_{1j} \ldot W_{mj}), 1\le j\le k$ of the random matrix $\mathcal{W}$
has stochastic representation
 $$(W_{1j} \ldot W_{mj}) = R_j\vk{O}_j=R_j(O_{1j},\cdots,O_{mj}),$$
 where $R_j$ and $\vk{O}_j$ are independent for any $ 1\le j\le k$. Here $\vk{O}_1 \ldot \vk{O}_k$ are independent copies of the random vector $\vk{O}$ which is uniformly distributed on the unit sphere of $\R^m$.
It follows that for any  $v$ in the support of $\zeta_1$ %if $\rho_i\not=0, i\le k$
\BQN\label{Saralees}
\Bigl(   \zeta_2 \ldot  \zeta_{k+1} \Bigr)  \Bigl \lvert (  \zeta_1=v) &\equaldis &
\Biggl( \sum_{i=1}^m( \rho_1 u_i + \widetilde{V}_1)^2  \ldot \sum_{i=1}^m( \rho_k u_i + \widetilde{V}_k)^2 \Biggr),
%\Bigl \lvert (  \zeta_1=v),
\EQN
where $\widetilde{V}_j= R_j O_{1j}, 1\le j\le k$ with  $u_j, 1\le j\le k$  are such that $\sum_{j=1}^m u_j^2=v$; the proof of \eqref{Saralees} is given in Appendix.
 A direct implication of \eqref{Saralees} is that %for any $v\not=0$ and $\rho_i\not=0, i\le k$
\BQNY\label{Saralees2}
\widetilde{\vk{\zeta}_v} \equaldis \Biggl( \frac{  \zeta_2 - \rho_1^2 v}{2 \rho_1 \sqrt{v}} \ldot \frac{ \zeta_{k+1}^2 - \rho_k^2 v}{2 \rho_k \sqrt{v}} \Biggr)  \Bigl \lvert (  \zeta_1=v) &\todis  &
\Bigl( \widetilde{V}_1 \ldot \widetilde{V}_k\Bigr), \quad v\to \IF. %+ O_p(1),
\EQNY
Consequently, \eqref{eq:Th1U1} holds with $\vk{U} \equaldis ( \widetilde{V}_1 \ldot \widetilde{V}_k)$.\\
Note that
if, for any $ 1\le j\le k$, $R_j$ is positive with $R_j^2$ having a chi-square distribution with $m$ degrees of freedom, then $\widetilde{V}_j$ is a $N(0,1)$ random variable.

Finally, we mention an extension of \netheo{Th2B}. It is possible therein to drop the assumptions that the rows of the matrix $\mathcal{W}_v$ have the same distribution. To this end, the condition \eqref{eqProb} needs to be re-stated, requiring the convergence of $w(v) \sqrt{v} \mathcal{W}_v$
to some random matrix $\mathcal{U}$. We shall illustrate the ideas in our first application below.

\section{Applications}
As it has been seen \tb{in} the Introduction,
conditional limit results are important in certain theoretical and applied models.
In this section, we shall present three applications of our main results.
The first \tb{one} concerns the derivation of Berman's sojourn limit theorems and the tail asymptotic behaviour of the supremum for a class of time-changed stationary chi-square processes.
In \tb{the} second application we shall investigate the maxima of perturbed
chi-square triangular arrays establishing both the convergence of the maxima
and a density type convergence result. Finally, motivated by the findings of
\cite{AsmRojas2008},
we consider the tail asymptotics of aggregated log-chi risks as our third  application.\\

\underline{\bf Berman's Sojourn Limit Theorem and Extremes of Time-Changed Chi-square Processes}:\\
Consider $\{X_i(t),t\ge 0\}, 1\le i\le m$ to be $m$ independent centered stationary Gaussian processes with %a.s. continuous sample paths and
covariance functions $r_i(\cdot), 1\le i\le m$ satisfying
\BQN\label{corrr}
 r_i(t)= 1 -  \cE{C_i}\abs{t}^{\alpha}+ o(\abs{t}^{\alpha}), \quad t\to 0, \quad \quad r_i(t)< 1, \quad \forall t> 0,
\EQN
with $\alpha \in (0,2]$ \cE{and $C_i, 1\le i \le m$ given positive constants.} Define  a time-changed stationary chi-square process $\{\zeta(t),t\ge 0\}$ by
$$ \zeta(t)= \sum_{i=1}^m X_i^2(\Theta_i t) , \quad t\ge 0,$$
where $\vk{\Theta}=(\Theta_1 \ldot \Theta_n)$ is a random vector with non-negative and bounded components being independent of the processes
$X_i,1\le i\le m$. We remark that time-changed processes are
\tb{used extensively}; see e.g., \cite{DHJ13, GemanYor01} and references therein. Next, let $\{Z_i(t),t\ge0\},1\le i\le m$ be independent copies of  a fractional Brownian motion $\{Z(t),t\ge 0\}$
with Hurst index $\alpha/2\in(0,1]$, i.e., $Z$ is a centered   Gaussian process with covariance function
$$ Cov(Z(s), Z(t))= t^\alpha+ s^\alpha- \abs{t-s}^\alpha, \quad s,t\ge 0.$$
We obtain below a conditional limit result which is crucial for derivation of Berman's sojourn limit theorems and the tail asymptotic behaviour of the supremum for the time-changed stationary chi-square processes.
Since $\zeta(0)$ has a chi-square distribution, it follows that its distribution function $G$ is in the Gumbel max-domain of attraction with scaling function $w(v)\equiv1/2$. We have, for any $0< t_1 < t_2 < \cdots< t_d$ and $x>0$ (set $\Delta_i(t_j)=
X_i(q(v)t_j)- r_i(q(v) t_j) X_i(0), X_{i,v}(t_j):=r_i(q(v) t_j) X_i(0) $ and $q(v)= v^{-1/\alpha}$)
\BQNY
\lefteqn{\Biggl( w(v) ( \zeta(q(v)t_1)- v)\ldot
 w(v) ( \zeta(q(v)t_d)- v)\Biggr) \Biggl\lvert ( \zeta(0)=v+ x/w(v))}\\
 &\equaldis & \Biggl( \frac{1}{2} \sum_{i=1}^m ( \Delta_i(\Theta_i t_1))^2
+  \sum_{i=1}^m  \Delta_i(\Theta_i t_1) X_{i,v}(\Theta_i t_1) +  \frac{1}{2} \Bigl( \sum_{i=1}^m (X_{i,v}(\Theta_i t_1))^2 - \zeta(0)\Bigr)  +x \ldot\\
&&
\frac{1}{2}  \sum_{i=1}^m ( \Delta_i(\Theta_i t_d))^2 + \sum_{i=1}^m  \Delta_i(\Theta_i t_d) X_{i,v}(\Theta_i t_d)
 +\frac{1}{2} \Bigl( \sum_{i=1}^m ( X_{i,v}(\Theta_i t_d))^2 - \zeta(0)\Bigr)+x \Biggr) \Biggl\lvert ( \zeta(0)=v+ 2x).
 \EQNY
By \eqref{corrr} it follows that
$$
 \Delta_i(\Theta_i t_j)\to 0 \ \text{and} \ r_i(q(v) \Theta_i t_j)\to 1 \ \ \text{almost\ surely} \ \text{as}\ u\to\IF,\ \ \ \forall 1\le i\le m,1\le j\le d.
$$
Consequently, the independence of $ \Delta_i(\Theta_i t_j)$ and $\zeta(0)$ implies
 \BQNY
 \lefteqn{\Biggl( w(v) ( \zeta(q(v)t_1)- v)\ldot
 w(v) ( \zeta(q(v)t_d)- v)\Biggr) \Biggl\lvert ( \zeta(0)=v+ x/w(v))}\\
 &\equaldis & \Biggl( O_p(1)
+  (1+o(1))\sum_{i=1}^m  \Delta_i(\Theta_i t_1) X_i(0)  -  \frac{1}{2} \sum_{i=1}^m (1-(r_i(q(v)\Theta_i t_1))^2)(X_i(0))^2  +x \ldot\\
&&
O_p(1)
+  (1+o(1))\sum_{i=1}^m  \Delta_i(\Theta_i t_d) X_i(0)  -  \frac{1}{2} \sum_{i=1}^m (1-(r_i(q(v)\Theta_i t_d))^2)X_i(0) +x \Biggr) \Biggl\lvert ( \zeta(0)=v+ 2x)\\
\EQNY
Furthermore, since $(X_1(0) \ldot X_m(0))$ is a standard Gaussian random vector, we have the stochastic representation
$$
(X_1(0) \ldot X_m(0))\equaldis R(O_1 \ldot O_m),
$$
where $(O_1 \ldot O_m)$ is a random vector uniformly distributed on the unit sphere of $\R^m$ being further independent of $R>0$ which is such that $R^2$ has a chi-square distribution with $m$ degrees of freedom. Here $R$ and $(O_1 \ldot O_m)$ are chosen such that they are independent of all the random variables (vectors) which are independent of $(X_1(0) \ldot X_m(0))$. Hence, in view of the independence between the random  variables (or vectors) we conclude that
\BQNY
 \lefteqn{\Biggl( w(v) ( \zeta(q(v)t_1)- v)\ldot
 w(v) ( \zeta(q(v)t_d)- v)\Biggr) \Biggl\lvert ( \zeta(0)=v+ x/w(v))}\\
 &\equaldis & \Biggl( O_p(1)
+  \sqrt{v+2x} \sum_{i=1}^m  \Delta_i(\Theta_i t_1) O_i  -  \frac{1}{2} (v+2x)\sum_{i=1}^m (1-(r_i(q(v)\Theta_i t_1))^2) O_i^2 +x \ldot\\
&&
O_p(1)
+ \sqrt{v+2x} \sum_{i=1}^m \Delta_i(\Theta_i t_d)O_i   -  \frac{1}{2} (v+2x)\sum_{i=1}^m (1-(r_i(q(v)\Theta_i t_d))^2)O_i^2 +x \Biggr) \Biggl\lvert ( R^2 =v+ 2x)\\
 &\todis & \Biggl(
 \sum_{i=1}^m  Z_i(C_i^{1/\alpha}\Theta_i t_1) O_i  -  \sum_{i=1}^m  C_i O_i^2\Theta^\alpha_i t_1^\alpha   +x \ldot
 \sum_{i=1}^m  Z_i(C_i^{1/\alpha}\Theta_i t_d) O_i  - \sum_{i=1}^m  C_i O_i^2\Theta^\alpha_i t_d^\alpha  +x \Biggr)
\EQNY
as $v\to \IF$, where $Z_i, O_i, 1\le i\le m, \vk{\Theta}$ are independent random elements.
Consequently, in view of \netheo{Th2B} we have the weak convergence of finite dimensional distributions
$$ \frac{1}{2}(\zeta(q(v)t) - v)  \Bigl\lvert (\zeta(0) > v) \todis
\widetilde{ Z(t)}:=
\sum_{i=1}^m  Z_i(C_i^{1/\alpha}\Theta_i t) O_i  -  \sum_{i=1}^m  C_i O_i^2\Theta^\alpha_i t^\alpha+ E, \quad t\ge 0,$$
where $E$ is a unit exponential random variable which is further independent of all the other random elements. Note that if $C_i=C\in (0,\IF), 1\le i\le m$ and $\vk{\Theta}$ has all components equal to 1, then
$$ \widetilde{Z(t)} \equaldis Z(C^{1/\alpha}t)- C t^{\alpha}+E, \quad t\ge 0,$$
which agrees with the findings of \cite{Berman82}.

Define the sojourn time of the process $\zeta$ above a level $v$ in the interval $[0,t]$ by
\BQN\label{eq:soj}
L_t(v)= \int_0^t \vk{1}(\zeta(s)> v) \, ds, \ \ t>0,
 \EQN
 where $\vk{1}(\cdot)$ is the indicator function. By checking the Assumptions in Theorem 3.1 in \cite{Berman82} (as it was done in Theorem 10.1 therein), we obtain the following  Berman's sojourn limit theorem for the time-changed stationary chi-square processes.

\BS Let $\{\zeta(t),t\ge0\}$ be the time-changed stationary chi-square process
with covariance functions satisfying \eqref{corrr},
and let $L_t(v)$ be defined as in \eqref{eq:soj}. Then, \tb{for all $t>0$ small,}
\BQN
\limit{v}  \int_x^\IF  \frac{ \pk{v^{1/\alpha}L_t(v)> y  }}{ v^{1/\alpha}\E{ L_t(v)}}\, dy =  B(x)
\EQN
holds at all continuity points $x>0$ of  $B(x)= \pk{\int_0^\IF \vk{1}(\widetilde{Z(s)}> 0)\, ds > x}$.
\ES

The tail asymptotics of supremum of  chi-square processes is important for
statistical analysis; see e.g.,  \cite{MR1043939} where
general stationary processes are considered.
Furthermore, the tail asymptotics of supremum of general (non-stationary)
chi-square processes are \tb{analyzed} in \cite{Pit96}. As mentioned in the Introduction, conditional limit results are crucial for establishing the asymptotic behaviour of supremum of stationary processes. Our approach for establishing the conditional limit result is direct, as shown above, and it differs from that of \cite{MR1043939}. Combined with the findings of \cite{MR1043939}, we have the following result.
\BS
Let $\{\zeta(t),t\ge0\}$ be the time-changed stationary chi-square process with
covariance functions satisfying \eqref{corrr}. Then, for any $T>0$,
\BQN
\pk{\sup_{t\in[0,T]}\zeta(t)>v}= H_\alpha[C_1 \ldot C_m] \frac{2^{1-m/2}T}{\Gamma(m/2)} v^{\frac{1}{\alpha}+\frac{m}{2}-1}\exp\left(-\frac{v}{2}\right)  (1+o(1)),
\EQN
as $v\to\IF$, where $\Gamma(\cdot)$ denotes the Euler Gamma function and
$$H_\alpha[C_1 \ldot C_m]=\lim_{a\downarrow 0}\frac{1}{a}\pk{\sup_{k\ge 1}\widetilde{Z(ak)}\le0} \in (0, \IF).$$

\ES

\underline{\bf Maxima of perturbed chi-square triangular arrays}:
We write below $H_\lambda$ for the  bivariate H\"usler-Reiss max-stable distribution,
\tb{defined as}
\BQNY
H_\lambda(x,y)&=& \exp\Biggl( - e^{-x}\Phi\Bigl(  \frac{\sqrt{\lambda}}{2} +\frac{y-x}{ \sqrt{\lambda}} \Bigr) -e^{-y}
\Phi\Bigl(  \frac{\sqrt{\lambda}}{2}+ \frac{x-y}{
\sqrt{\lambda}} \Bigr)  \Biggr), \quad x,y\inr,
\EQNY
with $\lambda \in (0,\IF)$ the dependence parameter, and $\Phi(\cdot)$ the standard Gaussian distribution function.
\tb{This distribution appeared initially in \cite{bro1977}, and was later studied in \cite{hue1989}.} It follows that the pdf $h_\lambda$ of $H_\lambda$ can be written as 
\BQN \label{hLL}
h_\lambda(x,y)= e^{-x}H_ \lambda(x,y) \Biggl(\frac{1 }{\sqrt{\lambda}}
\varphi\Bigl(  \frac{\sqrt{\lambda}}{2} +\frac{y-x}{ \sqrt{\lambda}} \Bigr) +
e^{ -y}\Phi\Bigl(  \frac{\sqrt{\lambda}}{2} +\frac{y-x}{ \sqrt{\lambda}}\Bigr)\Phi\Bigl(  \frac{\sqrt{\lambda}}{2}+ \frac{x-y}{
\sqrt{\lambda}} \Bigr)  \Biggr), \quad x,y\inr,
\EQN
with $\varphi(\cdot)$ the pdf of $\Phi(\cdot)$.\\
 The seminal contribution \cite{hue1989} shows that for bivariate Gaussian triangular arrays $(X_{i1}^{(n)},X_{i2}^{(n)}),
1 \le i\le n, n\ge 1$ where $(X_{i1}^{(n)},X_{i2}^{(n)}),1\le i\le n$ are independent bivariate Gaussian random vectors with $N(0,1)$ margins
and correlation $\rho_n\in (-1,1)/\{0\}$, the componentwise maxima is attracted by $H_\lambda$ if %the  H\"usler-Reiss
%condition
\BQN\label{eq:HRcondi}
 \limit{n} 4 \ln n (1- \rho_n)= \lambda\in[0,\IF).
 \EQN
  Let $H_n$ denote the joint distribution function
of a bivariate chi-square random vector $(\zeta_{1}^{(n)},\zeta_2^{(n)})$ as defined in \eqref{chi12}, where in \eqref{eJ} we put $\rho_n \in (-1,1)/\{0\}$ instead of $\rho$.
%with stochastic with margins that
%have chi-square distribution with $m$ degrees of freedom, and dependence parameter
%$\rho_n\in (-1,1)$ as in the Introduction.
In \cite{HashKab} the result of \cite{hue1989} was extended to chi-square case proving that under the condition \eqref{eq:HRcondi} (set below $t_n(x)= a_n x+ b_n$)
\BQN\label{convdf}
\limit{n} \sup_{x,y\inr} \Abs{ (H_n(t_n(x), t_n(y)))^n- H_\lambda(x,y)}=0,
\EQN
with
\BQN\label{anbN1}
a_n=2,\ \ b_n= 2 \ln n + (m-2) \ln (\ln n) - 2 \ln \Gamma(m/2).
\EQN
\COM{with $H_n$ the joint distribution function
of a bivariate chi-square random vector with margins that
have chi-square distribution with $m$ degrees of freedom, and dependence parameter
$\rho_n\in (-1,1)$.} Later on, in Theorem 2.2 in \cite{RevStat14} the same result for a perturbed chi-square vector were obtained, where the Gaussian assumption on $X_{i1}, 1\le i\le m$ in \eqref{chi12} is removed.  Instead therein both marginals $H_{n,j},i=1,2$  of $H_n$ are assumed to be  in the Gumbel max-domain of attraction, i.e.,
\BQN\label{saral}
\limit{n} \sup_{x\inr} \Abs{ (H_{n,j}(t_n(x)))^n- \exp(-\exp(-x))}=0,\quad j=1,2,
\EQN
where %for constants $a_n>0, b_n \inr$ given by
\BQN \label{anbN}
a_n= 1/ w(b_n)=2(1+o(1)), \quad b_n=G^{-1}(1- 1/n),
\EQN
with $G=H_{n,1}$ the distribution function of $\zeta_{1}^{(n)}$, %Note in passing that $H_{n,1}(x)=G(x), x\inr$.
\COM{Clearly, when $G$ is a chi-square distribution, we have
\BQN \label{w12}
w(x) \sim 1/2, \quad x\to \IF,
\EQN}
and further the H\"usler-Reiss condition
\BQN\label{r12}
 \limit{n}2 \frac{b_n}{a_n} (1- \rho_n^2)= \lambda \in [0,\IF)
 \EQN
 holds. Under the conditions \eqref{anbN} and  \eqref{r12},
  we have by \eqref{eq:Th2B} and Remarks \ref{RemTh2}, b) that %with $V$ a $N(0,1)$ random variable independent of $E$ being a unit exponential random variable
\BQNY
n \pk{ \zeta_{1}^{(n)} > t_n(x), \zeta_{2}^{(n)} > t_n(y)} &= &
\frac{\pk{\zeta_{1}^{(n)} > t_n(x)}}{\pk{\zeta_{1}^{(n)} > b_n}}  \pk{ \frac{ \zeta_{2}^{(n)}- t_n(x)}{a_n} >  y-x \Bigl \lvert \zeta_{1}^{(n)} > t_n(x) }(1+o(1)) \\
 &\to & \exp(-x)\pk{ \sqrt{\lambda} V  - \lambda/2 +E> y-x}, \quad n\to \IF,
\EQNY
where $V$ is a $N(0,1)$ random variable independent of the unit exponential random variable $E$.
This together with \eqref{saral} implies \eqref{convdf}, and thus the claim of Theorem 2.2 in \cite{RevStat14} follows. The result stated in \eqref{eq:Th2A} can be utilised to extend the convergence of distribution functions \eqref{convdf} to a convergence of the corresponding pdf's;
 see e.g., \cite{Faletal2010} for discussions on \tb{the convergence of densities}.

\BS \label{p:last}  Let $(\zeta_{1}^{(n)},\zeta_2^{(n)}), n\ge 1$ be a family of  bivariate chi-square random vectors  defined as in \eqref{chi12} with joint distribution function $H_n(x,y)$, where in \eqref{eJ} we put $\rho_n \in (-1,1)/\{0\}$ instead of $\rho$.
If  %\eqref{saral}, \eqref{anbN} and
\eqref{r12} is satisfied with $a_n$ and $b_n$ in \eqref{anbN1}
and $\hat{h}_n(x,y)$ is the pdf of $(H_n(t_n(x), t_n(y)))^n$ with $t_n(x)= a_n x+ b_n$, then
\COM{If,
 for all $x,y\inr$
 \BQN\label{26}
 \pk{  \zeta_{1}^{(n)}  > t_n(x) \Bigl \lvert  \zeta_{2}^{(n)}  = t_n(y)} = \pk{  \zeta_{2}^{(n)}  > t_n(y) \Bigl \lvert
\zeta_{1}^{(n)}  = t_n(x)}(1+o(1)), \quad n\to \IF,
\EQN
then we have}
\BQN
\limit{n} \hat{h}_n(x,y)& = &h_{\lambda}(x,y)
\EQN
holds for any $x,y\inr$.
\ES

%\begin{remark}
% The assumptions of \neprop{p:last} are satisfied for the case of chi-squared triangular arrays discussed in \cite{HashKab}.
%\end{remark}

\underline{\bf Aggregation of log-chi risks}:
Let $k\ge2,$ and define $\vk{ \zeta}:=(\zeta_1 \ldot \zeta_{k})$ to be a $k$-dimensional chi-square risk with $m$ degrees of freedom defined as in \eqref{eq:zetak} where $W_{i,j}$ in \eqref{eKJ} has variance $1-\rho_{j}^2\in(0,1)$ for any $1\le i\le m, 1\le j\le k-1$. Define further
$I_1 \ldot I_k$ to be independent of $\vk{ \zeta}$ and i.i.d.  Bernoulli random variables with
$\pk{I_i=1}=p=1- \pk{I_i=-1}$ and $p\in (0,1]$. For any constants $\sigma_j>0,   \mu_j\inr, 1\le j\le k$, define a $k$-dimensional log-chi risk $\vk{Z}=(Z_1 \ldot Z_k)$ by
\BQNY
Z_j= \exp(\sigma_j I_j \sqrt{\zeta_j}+\mu_j),\ \ \ 1\le j\le k.
\EQNY
The introduction of log-chi risks is motivated by  \cite{AsmRojas2008} where log-normal risks were considered, which are retrieved when  $m=1$ and $p=1/2$. As a generalization of the  result therein, we obtain the asymptotics of the aggregated log-chi risks.

\COM{
Consider $Z_1 \ldot Z_k, k\ge 2$ log-chi risks such that $Z_i= \exp(\sigma_i I_i \sqrt{\zeta_i}+\mu_i), i\le k$ where $ \zeta_i, i\le k$ are
given by \eqref{chim} being further independent of $I_i,i\le k$. The constants $\sigma_i>0,i\le k$ are assumed to be strictly positive, whereas
$\mu_i$'s are assumed to be finite. We shall suppose for simplicity that t
}
\BS\label{resmit} Let $Z_1 \ldot Z_k$ be log-chi risks  with $m$ degrees of freedom  as above. Let $\tilde \sigma:=\sigma_1 \ge \sigma_2 \ge \cdots \ge \sigma_k>0$, $\tilde \mu=\max_{1\le j\le k: \sigma_j=\tilde \sigma} \mu_j$, and
$
J_k=\sharp\{1\le j\le k: \sigma_j=\tilde \sigma, \mu_j=\tilde \mu\}
$. Then
\BQN
%\pk{\sum_{i=1}^n X_i> u} \sim J_k \pk{ \exp(\sigma_1  \zeta_1+ \mu)> u}= J_k \pk{  \zeta_1> (\ln u - \mu)/\sigma_1}, \quad u \to \IF,
\pk{\sum_{j=1}^k Z_j> u} &= & %pJ_k \frac{ 2^{m/2-1}}{\Gamma(m/2)} {t_u}^{m/2-1} \exp(- t_u^2/2), \quad u\to \IF,
 \frac{ pJ_k}{2^{m/2-1} \Gamma(m/2){{\tilde \sigma}}^{m- 2}} (\ln u - \tilde \mu)^{m- 2} \exp\left(- \frac{(\ln u - \tilde \mu)^2}{2{\tilde \sigma}^2}\right)(1+o(1)), \quad u\to \IF.
%\exp\Bigl(- \frac{(\ln u - \mu)^2}{2 \sigma_1^2} \Bigr), \quad u \to \IF,
\EQN
% where $J_k$ is the number of the elements of the index set $\{i\le k: \sigma_i=\tilde \sigma, \mu_i=\tilde \mu\}$, the integer $m$ is
 %the multiplicity of $\tilde \sigma$ among $\sigma_1 \ldot \sigma_k$ and $t_u:= (\ln u - \tilde \mu)/\tilde \sigma$.
 \ES
Note in passing that %if $m=1$ and $p=1/2$ the above asymptotics agrees with that derived in  \cite{AsmRojas2008} for log-normal risks.
the tail asymptotics of the maximum $\max_{1 \le j\le k} Z_j$
can be further shown to be tail-equivalent with the total risk $\sum_{j=1}^k Z_j$; see \cite{Foss2011} for more examples on this topic.\\

\section{Proofs}

\prooftheo{Th2}
For any $v>0$ we have
\BQNY
\lefteqn{\Biggl( \frac{  \zeta_2 - \rho_1^2v}{2\rho_1\sqrt{v}} \ldot
\frac{ \zeta_{k+1} - \rho_k^2v}{2\rho_k\sqrt{v}}\Biggr) \biggl \lvert   (\zeta_1 =v)}\\
&\equaldis &
\Biggl( \frac{\sum_{i=1}^m (\rho_1 X_{i1} + W_{i1})^2  - \rho_1^2v}{2\rho_1\sqrt{v}} \ldot
\frac{\sum_{i=1}^m (\rho_k  X_{i1} + W_{ik})^2  - \rho_k^2v}{2\rho_k\sqrt{v}}\Biggr)\biggl \lvert (  \zeta_1 =v)\\
&\equaldis &
\Biggl( \frac{ 2\rho_1 \sum_{i=1}^m X_{i1} W_{i1} + \sum_{i=1}^m W_{i1}^2  }{2\rho_1\sqrt{v}} \ldot
\frac{2\rho_k \sum_{i=1}^m X_{i1} W_{ik} + \sum_{i=1}^m W_{ik}^2}{2\rho_k\sqrt{v}}\Biggr)\biggl \lvert (  \zeta_1 =v).
%\\
%&\equaldis &
%\Biggl( \frac{ 2\rho_1 \sqrt{ \sum_{i=1}^m X_{i1}^2} W_{i1} + \sum_{i=1}^m W_{i1}^2  }{2\rho_1\sqrt{v}} \ldot
%\frac{2\rho_k \sum_{i=1}^m X_{i1} W_{ik} + \sum_{i=1}^m W_{ik}^2}{2\rho_1\sqrt{v}}\Biggr)\biggl \lvert (  \zeta_1 =v)\\
\EQNY
Since further by the independence of $(X_{11} \ldot X_{m1})$ and the Gaussian random matrix $\mathcal{W }$ we have
\BQN\label{babo}
\Biggl( \sum_{i=1}^m X_{i1} W_{i1}  \ldot
\sum_{i=1}^m X_{i1} W_{ik} \Biggr)\biggl \lvert (  \zeta_1 =v) &\equaldis &
\Biggl( \sqrt{ \sum_{i=1}^m X_{i1}^2} W_{1} \ldot
\sqrt{\sum_{i=1}^m X_{i1}^2} W_{k} \Biggr)\biggl \lvert (  \zeta_1 =v) \notag\\
& \equaldis & (\sqrt{v} W_1 \ldot \sqrt{v} W_k).
\EQN
Thus, the first claim follows immediately by the fact that $\sum_{i=1}^m W_{ij}^2/\sqrt{v} \toprob 0$ for any $1\le j\le k$.\\ %almost surely as $v\to \IF$
Next, the assumption that  $G$ of $  \zeta_1$ is in the Gumbel max-domain of attraction,
implies  $\limit{v} v w(v)=\IF$ and the convergence in distribution
$$ w(v)(  \zeta_1- v) \lvert (  \zeta_1 > v) \todis E, \quad v\to \IF.$$
By \tb{the above} we obtain (set $v_z:= v+ z/w(v)$) %as $v\to \IF$
\BQNY
\lefteqn{\pb{ \frac{  \zeta_2 - \rho_1^2v}{2\rho_1\sqrt{v}} \le x_1    \ldot
\frac{ \zeta_{k+1} - \rho_k^2v}{2\rho_k\sqrt{v}} \le x_k \biggl \lvert    \zeta_1= v+ z/w(v)}}\\
 &= &
\pb{ \frac{   \zeta_2 - \rho_1^2v_z + \rho_1^2 z /w(v) } { 2\rho_1 \sqrt{v_z}} \sqrt{v_z/v} \le x_1    \ldot
\frac{ \zeta_{k+1} - \rho_k^2 v_z+ \rho_k^2 z/w(v)}{ 2\rho_k\sqrt{v_z}} \sqrt{v_z/v} \le x_k \biggl \lvert    \zeta_1= v_z} \\
& = &\pb{ \frac{   \zeta_2 - \rho_1^2v_z } { 2\rho_1 \sqrt{v_z}}  +  \rho_1 c z\le x_1    \ldot
\frac{ \zeta_{k+1} - \rho_k^2 v_z}{ 2\rho_k\sqrt{v_z}}   +\rho_k c z \le x_k \biggl \lvert    \zeta_1= v_z}(1+o(1))\\
&\to& \pk{W_1 \le x_1-\rho_1 c z \ldot W_k \le x_k-\rho_kc z } , \quad v\to \IF,
\EQNY
where the convergence  holds uniformly with respect to $z\inr$, and thus we can substitute $z$ by $z_v, v>0$ satisfying $\lim_{v\to \IF} z_v= z\inr$ in the above.
Consequently, in view of Lemma  4.2 in \cite{MR2301753}, we obtain
\BQN\label{eq:GG}
\lefteqn{\pb{ \frac{  \zeta_2 - \rho_1^2v}{2\rho_1\sqrt{v}} \le x_1    \ldot
\frac{ \zeta_{k+1} - \rho_k^2v}{2\rho_k\sqrt{v}} \le x_k \biggl \lvert    \zeta_1> v}}\nonumber\\
&=&\int_v^\IF \pb{ \frac{  \zeta_2 - \rho_1^2v}{2\rho_1\sqrt{v}} \le x_1    \ldot
\frac{ \zeta_{k+1} - \rho_k^2v}{2\rho_k\sqrt{v}} \le x_k \biggl \lvert    \zeta_1= s}\, d G(s)/ (1- G(v))\nonumber\\
&=& \int_{0}^\IF \pb{ \frac{  \zeta_2 - \rho_1^2v}{2\rho_1\sqrt{v}} \le x_1    \ldot
\frac{ \zeta_{k+1} - \rho_k^2v}{2\rho_k\sqrt{v}} \le x_k \biggl \lvert    \zeta_1= v+ z/w(v)} \, d G(v+ z/w(v))/(1- G(v))\nonumber\\
&\to  & \int_{0}^\IF  \pk{W_1 \le x_1-\rho_1 c z \ldot W_k \le x_k-\rho_kc z }  \exp(-z)\,  dz,
\quad v\to \IF\nonumber\\
&=& \pk{U_1 + \rho_1 c E \le x_1 \ldot U_k +\rho_kc E \le x_k },
\EQN
thus the proof is complete. \QED

\prooftheo{Th2B} First note that \eqref{vww} implies, as $v\to \IF$,
\BQNY
w(v) (v- \rho_{j,v}^2 v)= 2vw(v) (1- \rho_{j,v} )(1+o(1))  = \lambda_j/2, \quad 1\le j\le k.
\EQNY
Further, the scaling function $w$ is self-neglecting, i.e.,
$$\lim_{v\to\IF}\frac{w( v+ x/w(v))}{w(v)} = 1,\ \ \ \ \forall x\inr.$$ % uniformly on compact sets that contain $x$, then
Therefore, the claim of statement i) follows by the assumption \eqref{eqProb} and the convergence in distribution
\BQNY
\Biggl(  w(v)(\zeta_{2v}- \rho_{1,v}^2 v) \ldot   w(v)(\zeta_{k+1, v}- \rho_{k,v}^2v)\Biggr) \Bigl \lvert   (\zeta_{1}= v) \todis
( 2 U_1  \ldot 2 U_k )
\EQNY
which can \tb{be} confirmed as in \eqref{babo},
with the aid of the assumption \eqref{eqProb}.
The claim of statement ii) can be established using similar
arguments as in \eqref{eq:GG}.
\tb{This completes the proof.} \QED

\proofprop{p:last}
%We first introduce some further notation.
Denote by $h_n(x,y)$  the  pdf of $H_n(x ,y)$ and write $h_{n,j},j=1,2$ for the marginal pdf's of it. Further, denote $h_n( \cdot \lvert x)$ to be the conditional pdf of $\zeta_2^{(n)} \lvert \zeta_1^{(n)}= x$. %(see e.g., \cite{MR2830509})
By \netheo{Th2B} and Remarks \ref{RemTh2}, b), we have, for any $x,y\inr$
\BQNY
\limit{n} \pk{ \zeta_{2}^{(n)}  \le t_n(y)  \lvert \zeta_{1}^{(n)}=t_n(x)}= \pk{ V  \le  \frac{\sqrt{\lambda}}{2} +\frac{y-x}{\sqrt{\lambda} }} ,
\EQNY
with $V$ a $N(0,1)$ random variable. By symmetry
\BQNY
\limit{n} \pk{ \zeta_1^{(n)}  \le t_n(x)  \lvert  \zeta_2^{(n)}=t_n(y)}= \pk{ V  \le  \frac{\sqrt{\lambda}}{2} +\frac{x-y}{\sqrt{\lambda} }}.
\EQNY
Consequently, since $\limit{n} (H_{n}(t_n(x),t_n(y)))^{n}=H_\lambda(x,y)$
\BQNY
 \hat{h}_n(x,y)&=& a_n^2n(H_{n}(t_n(x),t_n(y)))^{n-1} h_{n}(t_n(y) \lvert t_n(x)) h_{n,1}(t_n(x))\\
&&+ a_n^2n(n-1)(H_{n}(t_n(x),t_n(y)))^{n-2}h_{n,1}(t_n(x))h_{n,2}(t_n(y))\\
&& \times  \pk{ \zeta_2^{(n)}  \le t_n(y)  \lvert  \zeta_1^{(n)}=t_n(x)}
\pk{ \zeta_1^{(n)} \le t_n(x)  \lvert  \zeta_2^{(n)}=t_n(y)}\\
&=&(1+o(1)) H_\lambda(x,y) \Biggl[ a_n  h_{n}(t_n(y)\lvert t_n(x)) a_n n h_{n,1}(t_n(x))\\
&& + e^{-x+y} \pk{ V  \le  \frac{\sqrt{\lambda}}{2} +\frac{y-x}{\sqrt{\lambda} }} \pk{ V  \le \frac{\sqrt{\lambda}}{2} +\frac{x-y}{\sqrt{\lambda} }}  \Biggr].
\EQNY
Since $G$ is a chi-square distribution, we have  $n a_n h_{n1}(t_n(x)) \to \exp(- x)$ as $n\to \IF$.
%Next, note that (18) holds also in the setup of this proposition, and by that result the conditional distribution is the same as in the case of chi-square triangular arrays, consequently,
Further, in the light of \tb{Remarks \ref{RemTh2}, c)} we obtain that
$$ \limit{n} a_n h_{n}(t_n(y) \lvert t_n(x))=  g(y \lvert x), $$
with $g( \cdot \lvert x)$ the pdf of $\sqrt{\lambda} V - \lambda/2 +x$, implying thus
\BQNY
\limit{n}\hat{h}_n(x,y)&=& e^{-x}H_\lambda(x,y) \Biggl[ \frac{1}{ \sqrt{ \lambda} } \varphi( \frac{\sqrt{\lambda}}{2} +\frac{y -x}{\sqrt{\lambda} }  )
+ e^{-y} \Phi(\frac{\sqrt{\lambda}}{2} +\frac{y-x}{\sqrt{\lambda} })
 \Phi(\frac{\sqrt{\lambda}}{2} +\frac{x-y}{\sqrt{\lambda} })  \Biggr]
=h_\lambda(x,y),
\EQNY
hence the proof is complete.\QED

\def\tilsig{\tilde \sigma}
\def\tilmu{\tilde \mu}
\def\tilZ{\tilde Z}

\proofprop{resmit} %It followsFrom \cite{MR3003975}
The proof is based on Theorem 4.2 in \cite{MitraResn2009}. Let $\tilZ= \exp(\tilsig I_1 \sqrt{\zeta_1}+\tilmu)$.
Since $\sqrt{\zeta_1}$ is a chi-distribution with $m$ degrees of freedom, it follows that
$$
\pk{\tilde Z>u}=\frac{ p }{2^{m/2-1} \Gamma(m/2){{\tilde \sigma}}^{m- 2}} (\ln u - \tilde \mu)^{m- 2} \exp\left(- \frac{(\ln u - \tilde \mu)^2}{2{\tilde \sigma}^2}\right)(1+o(1)), \quad u\to \IF,
$$
implying that $\tilZ$ has distribution function in the Gumbel max-domain of attraction with scaling function $w(x)= \ln x/(\tilsig^2 x)$.
Since further $\lim_{x\to \IF} w(x)= 0$, in view of Theorem 4.2 in \cite{MitraResn2009} we conclude
the claim  by checking Assumptions 2.3-2.5 therein. In our setup it suffices to show them for $k=2$.
\tb{ For the simplicity  of presentation, we assume further that} $\sigma_1=\sigma_2=1$, $p=1$ and $\mu_1=\mu_2=0$. For any $a>0$ we have
\BQNY
\pk{ w(u) Z_2> a \lvert Z_1 > u}&=&  \pk{  Z_2> a u/\ln u  \Bigl\lvert Z_1 > u}\\
%&=&  \pk{   \zeta_2 > \ln (a u/\ln u)  \lvert  \zeta_1 > \ln u}\\
%&=&  \pk{   \zeta_2 > \ln a + \ln u - \ln \ln u)  \lvert  \zeta_1 > \ln u}\\
%&=&  \pk{   \zeta_2 > \ln a + \sqrt{v} - \ln \sqrt{v})  \lvert  \zeta_1 > \sqrt{v}}\\
&\overset{v=(\ln u)^2}=&  \pk{    \zeta_2 > (\ln a + \sqrt{v} - \ln \sqrt{v})^2   \Bigl\lvert   \zeta_1 > v}\\
%&=&  \pk{    \zeta_2 > v \Bigl(1  + (\ln a   - \ln \sqrt{v})/\sqrt{v}\Bigr) ^2   \Bigl\lvert   \zeta_1 > v}\\
&=&  \pk{  \frac{  \zeta_2- \rho_1^2 v}{\sqrt{v}} > \sqrt{v} \Bigl[ \Bigl(1  + (\ln a   - \ln \sqrt{v})/\sqrt{v}\Bigr) ^2 - \rho_1^2\Bigr]   \Bigl\lvert \zeta_1 > v}
%\\
\to  0%, \quad v\to \IF,
\EQNY
 as $u\to \IF,$ where the last convergence follows from \netheo{Th2}  and the fact that $\rho_1^2< 1$, hence Assumption 2.3
 and Assumption 2.4 (by symmetry) in \cite{MitraResn2009} hold. The Assumption 2.5 in  \cite{MitraResn2009} follows if we show that
\BQNY
\frac{\pk{ \min (Z_1, Z_2)> u/\ln u}}{ \pk{ Z_1> u}} =
\frac{\pk{ \min (  \zeta_1,   \zeta_2)> v^*}}{ \pk{   \zeta_1> (\ln u)^2}}
&\to & 0%, \quad u\to \IF,
\EQNY
as $u \to \IF$, where $v^*=(\ln u - \ln \ln u)^2=(\ln u)^2(1+o(1))$.
By the definition of $( \zeta_1, \zeta_2)$ we have the stochastic representation
$$  \zeta_1+   \zeta_2\equaldis  (1+ \rho_1) \sum_{i=1}^m W_i^2+  (1-\rho_1) \sum_{i=m+1}^{2m} W_i^2, $$
where $W_i,1\le i\le 2m$ are i.i.d. $N(0,1)$ random variables. Without loss of generality we may assume that $\rho_1>0$.
Let $c:=2/(1+\rho_1)>1$. We have
\BQNY
\pk{ \min (  \zeta_1,   \zeta_2)> u} &\le& \pk{   \zeta_1+  \zeta_2> c(1+\rho_1)u }\\
&=&(1+ 1/\rho_1)^{m/2}\frac{2^{1- m}c^{m/2-1}}{\Gamma(m/2)} u^{m/2-1} \exp\left(- \frac{c u}{2}\right) (1+o(1))
\EQNY
as $ u\to \IF$, where the last asymptotic equivalence follows from
 Example 5 in \cite{HashKorshPit}. Consequently,
\BQNY
\limit{u}\frac{\pk{ \min (Z_1, Z_2)> u/\ln u}}{ \pk{ Z_1> u}}\le\limit{v}\frac{\pk{   \zeta_1+  \zeta_2> c(1+\rho_1)v^* }}{ \pk{   \zeta_1> v^*(1+o(1))}}=0
\EQNY
and thus the proof is complete.  \QED

\def\tHH{\rho_*}
\section{Appendix}
In this section we first present the proof of \eqref{Saralees} and then give a direct proof of \eqref{eq:A}.\\
For notational simplicity we consider only the case $k=2$. We have, for any $v$ in the support of $\zeta_1$
\BQNY
\Bigl(   \zeta_2,  \zeta_3\Bigr) \Bigl\lvert (  \zeta_1=v)
& \equaldis &
\Bigl( \sum_{i=1}^m (\rho_1 X_{i1}+ W_{i1})^2,  \sum_{i=1}^m (\rho_2 X_{i1}+ W_{i2})^2\Bigr) \Biggl\lvert (  \zeta_1=v)\\
& \equaldis &
\Biggl( \rho_1 ^2 v + 2 \rho_1 \sum_{i=1}^m X_{i1}W_{i1} + \sum_{i=1}^m W_{i1} ^2,
\rho_2^2 v+ 2 \rho_2 \sum_{i=1}^m  X_{i1} W_{i2}+ \sum_{i=1}^m W_{i2}^2\Biggr)\Biggl \lvert (  \zeta_1=v).
\EQNY
The assumption that %Since $(W_{11} \ldot W_{1m}) $ and $(W_{21} \ldot W_{2m})$ are independent with
$(W_{1j} \ldot W_{mj}) \equaldis R_j \vk{O}_j, j=1,2$ %with $\vk{O}_1, \vk{O}_2$ having the uniform distribution on the unit sphere of $\R^m$ and being are independent
implies for any $v$ in the support of $\zeta_1$
\BQNY
\Bigl(   \zeta_2,  \zeta_3\Bigr) \Bigl\lvert (  \zeta_1=v)
& \equaldis &
\Biggl( \rho_1 ^2 v + 2 \rho_1 R_1 O_{11}\sqrt{\sum_{i=1}^m X_{i1}^2}  + R_1^2,
\rho_2^2 v+ 2 \rho_2 R_2 O_{12} \sqrt{\sum_{i=1}^m  X_{i1}^2} + R_2^2\Biggr) \Biggl\lvert \Bigl(  \sum_{i=1}^m  X_{i1}^2=v\Bigr)\\
& \equaldis &
\Biggl( \rho_1 ^2 v + 2 \rho_1 R_1 O_{11} \sqrt{v} + R_1^2,
\rho_2^2 v+ 2 \rho_2 R_2 O_{21} \sqrt{v} + R_2^2\Biggr) \\
& \equaldis &
\Bigl( \sum_{i=1}^m (\rho_1 u_i+ W_{i1})^2,  \sum_{i=1}^m (\rho_2 u_i+ W_{i2})^2\Bigr) ,
\EQNY
for any $u_i,i\le m$ such that $\sum_{i=1}^m u_i^2=v$, hence the claim  of \eqref{Saralees} follows.\\

Next, we show the proof of \eqref{eq:A}. In view of \cite{pre06111634}  (see also \cite{MR2674417}),  we have the stochastic representation
%\BQN\label{repA}
$$(  \zeta_1,   \zeta_2)\equaldis (U_m,V_m) ,$$
%\EQN
where (set $N=m+1$)
$$U_{m} = \sum_{i=1}^N \Bigl(X_{i1} - \overline{X_{N1}}\Bigr)^2, \quad V_{m}= \sum_{i=1}^N \Bigl(X_{i 2} - \overline{X_{N2}}\Bigr)^2, \quad
\overline{X_{N1}}:=\sum_{i=1}^N \frac{X_{i1}}{N}, \quad \overline{X_{N2}}:=\sum_{i=1}^N \frac{X_{i2}}{N},$$
which follows %The claim in \eqref{repA} follows
from the facts that
$(U_m, V_m)$ is independent of $(\overline{X_{N1}},\overline{X_{N2}})$, and
$(\overline{X_{N1}} \sqrt{N} , \overline{X_{N2}} \sqrt{N} )$ has the same distribution as $(X_{11},X_{12})$.
From equation (3) in \cite{MR2674417} %, distributional properties of $(  \zeta_1,   \zeta_2)$ can be easily studied since
we have the following expression for the pdf $h_m$ of $(U_{m},V_{m})$:
\begin{eqnarray}
\label{eq:cond}
h_m(u, v)= \frac {(uv)^{m/2-1}}{2^m (\Gamma(m/2))^2 \left( \tHH^2\right)^{m/2}}\exp\Bigl( -\frac {u+v}{2 \tHH^2} \Bigr)
\  {}_0F_1 \left( ; \frac{m}{2};  \frac{\rho^2 uv}{\left(2 \tHH^2\right)^2} \right),
\quad
\forall u, v \in (0,\IF),
\end{eqnarray}
where $\tHH:= \sqrt{1 - \rho^2}$ and $_0F_1(;a;x)= \sum_{n=0}^\IF \frac{\Gamma(a+n)}{\Gamma(a)} \frac{x^n}{n!}$.
 %with $\Gamma(\cdot)$ the Euler Gamma function.
  By \eqref{eq:cond} %, if $h_m(\cdot \lvert v), v>0$ denotes the pdf  of $  \zeta_2 \lvert (  \zeta_1=v)$,
%then
 the pdf $g_m(x \lvert v),x\inr$ of the conditional random variable %(set below $$)
 $$ \zmvv =\frac{  \zeta_2  - \rho^2v }{ \tHH  \sqrt{v}} \Bigl \lvert (  \zeta_1= v)$$
  is  given by (set $x_\rho:= \tilde\rho  x \sqrt {v} + \rho^2 v $ where $\tilde\rho:=1 - \rho^2$)
\BQNY
g_m(x\lvert v) %&=& \tilde\rho   \sqrt {v} h_m(  \tilde\rho  x \sqrt {v} + \rho^2 v  \lvert v)\\
&=&
\frac{ \tilde\rho   \sqrt {v} }{ \Gamma(m/2)( 2 \tilde\rho )^{m/2}} x_\rho^{(m-2)/2}
\exp\Biggl(- \frac{ x \sqrt {v}}{2 } \Biggr) \exp\Biggl(- \frac{ \rho^2 v }{\tilde\rho } \Biggr)
 \ _0F_1\Biggl(;\frac{m}{2}; \frac{  \rho^2  v x_\rho }{ (2 \tilde\rho)^2 }    \Biggr).
\EQNY
Utilising the well-known asymptotic expansion
\BQNY\label{exp0F1}
 _0F_1(;m;z) =\frac{\Gamma(m)}{ 2\sqrt{ \pi}} z^{1/4- m/2} \exp(2 \sqrt{z}) \Bigl(1+ O(1/\sqrt{z})\Bigr), \quad z\to \IF
 \EQNY
we can further write, as $v\to \IF$,
\BQNY
g_m(x\lvert v)
 &=&
 \frac{\tilde\rho   \sqrt {v} }{ \Gamma(m/2)( 2 \tilde\rho )^{m/2}} x_\rho^{(m-2)/2}
\exp\Biggl(- \frac{   x \sqrt {v}}{2  } \Biggr) \exp\Biggl(- \frac{ \rho^2 v }{2\tilde\rho} \Biggr)\\
&& \times
\frac{\Gamma(m/2)}{ 2\sqrt{ \pi}} \Biggl( \frac{  \rho^2  v x_\rho }{(2 \tilde\rho )^2} \Biggr)^{(1- m)/4}\exp
\Biggl(   \Biggl(\frac{  \rho^2  v x_\rho}{\tilde\rho ^2} \Biggr)^{1/2}\Biggr)(1+o(1))\\
% &=&
% \frac{1}{ 2\sqrt{ \pi}}
%\frac{ \sqrt {v} \tilde\rho  }{ \sqrt{  2 \tilde\rho  \rho v}}
%\exp\Biggl(- \frac{  x \sqrt {v}}{2 } \Biggr) \exp\Biggl(- \frac{ \rho^2 v }{ \tilde\rho  } \Biggr)
%\exp\Biggl( \frac{\rho^2 v}{ \tilde\rho  }  \Biggl(1+  \frac{ \tilde\rho x }{\rho^2\sqrt{v}} \Biggr)^{1/2}\Biggr)\\
% &\sim&
%\frac{1}{ 2\sqrt{ \pi}}
%\frac{  \sqrt{\tilde\rho}   }{ \sqrt{2 \rho}}   \exp\Biggl(- \frac{  x \sqrt {v}}{2  } \Biggr) \exp\Biggl(- \frac{ \rho^2 v }{ \tilde\rho  } \Biggr)\exp\Biggl( %\frac{\rho^2 v}{ \tilde\rho  }  \Biggl(1+   \frac{\tilde\rho x }{2 \rho^2 \sqrt{v}}  - \frac{1}{8}  \frac{ \tilde\rho ^2 x^2}{\rho^4 v}+ o(v)\Biggr)\Biggr)\\
% &\sim&
%\frac{1}{ \sqrt{ 2 \pi}}
%\frac{  \sqrt{\tilde\rho}   }{ 2 \sqrt{ \rho}}    \exp\Biggl( -\frac{  \tilde\rho x ^2 }{8 \rho^2}+ o(v)\Biggr)\\
 &=&\frac{1}{ \sqrt{ 2 \pi}}
\frac{ \sqrt{\tilde\rho}   }{ 2 \sqrt{ \rho}}    \exp\Biggl( -\frac{ \tilde\rho x ^2 }{8 \rho^2}\Biggr)(1+o(1)).
\EQNY
\def\tilde\rho{\widetilde{\rho}}
Consequently, for any $x\inr$
$$ \pb{  \frac{  \zeta_2  - \rho^2v  }{2 \rho  \sqrt{1- \rho^2}  \sqrt{v}} \le x \Bigl \lvert   \zeta_1= v}=  \pb{ \zmvv  \le \frac{2\rho x}{\sqrt{1- \rho^2}} }
\to \pk{W_1\le x}, \quad v\to \IF. $$\\

\textbf{Acknowledgments.} Partial support from the Swiss National Science Foundation Project 200021-1401633/1  and by the project RARE -318984
 (a Marie Curie International Research Staff Exchange Scheme Fellowship within the 7th European Community Framework Programme) is kindly acknowledged. The first author also acknowledeges partial support by NCN Grant No 2011/01/B/ST1/01521 (2011-2013).

\COM{
\newpage

Consider a Gaussian triangular array with correlation $\rho_n$ as in H\"usler-Reiss and suppose that

$$ 2(1- \rho_n) \ln n  \sim  (1- \rho_n) b_n/a_n  \sim  2 \lambda^2\in (0,\IF)$$
Let us find the limit of
$$ \pk{X_{1n} > a_n x+ b_n \lvert X_{2n}= a_n y+ b_n}$$
we have
$$ X_{1n}= \rho_n X_{2n}+ \sqrt{1 - \rho_n^2} Z$$
hence
\BQNY
\pk{ X_{1n} > a_n x+ b_n \lvert X_{2n}= a_n y+ b_n}&= &
\pk{ \rho_n (a_n y+ b_n) + \sqrt{1- \rho_n^2} Z > a_n x+ b_n }\\
&= &
\pk{ \frac{1}{a_n} \sqrt{1- \rho_n^2} Z >  x- \rho_n y + b_n(1- \rho_n)/a_n }\\
&= &
\pk{ \sqrt{2(1- \rho_n) b_n /a_n} Z >  x- \rho_n y + b_n(1- \rho_n)/a_n }\\
&\to &
\pk{ 2 \lambda Z >  x- y + 2 \lambda^2}\\
\EQNY
This limit yields the H\"usler-Reiss distribution
\BQNY
H_\lambda^*(x,y)&=& \exp\Biggl( - e^{-x}\Phi\Bigl(  \lambda+\frac{y-x}{ 2 \lambda} \Bigr) -e^{-y}
\Phi\Bigl(  \lambda+ \frac{x-y}{
2\lambda} \Bigr)  \Biggr), \quad x,y\inr,
\EQNY
If we put $2\lambda^2 = \tau/2$, or $\lambda= \sqrt{\tau}/2$ the above limit can be re-written as
\BQN\label{A}
\pk{ X_{1n} > a_n x+ b_n \lvert X_{2n}= a_n y+ b_n}
&\to &
\pk{ \sqrt{\tau}  Z >  x- y + \tau/2}
\EQN
and the limit distribution is
\BQNY
H_\tau(x,y)&=& \exp\Biggl( - e^{-x}\Phi\Bigl(  \frac{\sqrt{\tau}}{2} +\frac{y-x}{ \sqrt{\tau}} \Bigr) -e^{-y}
\Phi\Bigl(  \frac{\sqrt{\tau}}{2}+ \frac{x-y}{
\sqrt{\tau}} \Bigr)  \Biggr), \quad x,y\inr.
\EQNY

Now we need to prove that if $(X_{1n},X_{2n})$ is the chi-square triangular array that satisfies condition
$$ (1- \rho_n)^2 \frac{b_n}{a_n} \to \tau/2 $$
then $H_\tau$ is the limit distribution of the maxima. This follows if we show that \eqref{A} holds.

We have as in the proof of the first theorem (here $a_n=2$)

\BQNY
\pk{ X_{1n} > a_n x+ b_n \lvert X_{2n}= a_n y+ b_n}
&\sim &
\pk{ \rho_n^2 (a_n y+ b_n) + 2 \rho_n  \sqrt{a_n y+ b_n}  \sqrt{1- \rho_n^ 2} Z > a_n x+ b_n }\\
&\sim &
\pk{ 2 \sqrt{b_n}(1+o(1))  \sqrt{1- \rho_n^ 2} Z > 2 x- 2 y+ b_n(1- \rho_n^2) }\\
&\sim &
\pk{  \sqrt{b_n}(1+o(1))  \sqrt{1- \rho_n^ 2} Z >  x-  y+ b_n(1- \rho_n^2)/a_n }\\
&\sim & \pk{ \sqrt{ \tau} Z >  x- y +  \tau/2}\\
\EQNY
}

\bibliographystyle{plain}
\bibliography{chisq}
\end{document}